\documentclass[10pt]{article}
\usepackage{amsfonts}

\usepackage{amscd} 
\usepackage{graphicx}
\usepackage{amsmath}
\usepackage{amsfonts}
\usepackage{amssymb}

\newenvironment{proof}[1][Proof]{\textbf{#1.} }{\ \rule{0.5em}{0.5em}}

\renewcommand{\Im}{\operatorname{Im}}

\newtheorem{theo}{Th\'{e}or\`{e}me}[section]
\newtheorem{defi}{D\'{e}finition}[section]
\newtheorem{propo}{Proposition}[section]
\newtheorem{lemm}{Lemme}[section]

\newtheorem{coro}{Corollaire}[section]
\newtheorem{rmk}{Remarque}[section]

\newtheorem{exemple}{Exemple}[section]

\begin{document}

\title{Formes d'inertie et complexe de Koszul associ\'{e}s \`{a} des polyn\^{o}mes
plurihomog\`{e}nes\thanks{%
Ce travail a \'{e}t\'{e} \'{e}labor\'{e} avec l'aide de la coop\'{e}ration
franco-marocaine Action in\'{e}tgr\'{e}e A.I. MA/02/32 et du programm PAS
27/2001 financ\'{e} par l'AUPELF.} }
\author{Azzouz AWANE$^{\left( 1\right) }$, Abdelouahab CHKIRIBA$^{\left( 1\right) }$
et Michel GOZE$^{\left( \small 2\right) }$ }
\date{}
\maketitle
\begin{center}
{\small(1)}{\scriptsize UFR de G\'{e}om\'{e}trie Diff\'{e}rentielle et Applications%
}\\
{\scriptsize Facult\'{e} des Sciences Ben M'sik. B.P. 7955. Boulevard Driss
Harti. Casablanca. Maroc}\\
{\small(2)}{\scriptsize Facult\'{e} des Sciences et Techniques. Universit\'{e} de
Haute Alsace}\\
{\scriptsize 4, rue des Fr\`{e}res Lumi\`{e}re. F. 68093 Mulhouse Cedex }\\
{\scriptsize E-mail : awane@eudoramail.com . M.Goze@uha.fr}
\end{center}
\begin{abstract}
{\scriptsize The existence of common zero of a family of polynomials has led
to the study of inertial forms, whose homogeneous part of degree 0
constitutes the ideal resultant. The Kozsul and \v {C}ech cohomologies
groups play a fundamental role in this study. An analogueous of Hurwitz
theorem is given, and also, one finds a N.H.Mcoy theorem in a particular
case of this study. }
\end{abstract}



\section{Introduction}

La notion des formes d'inertie (Tr\"{a}gheitsformen) est due \`{a} F.Mertens 
\cite{MRTNS} au XIX-i\`{e}me$\ $si\`{e}cle, est li\'{e}e au probl\`{e}me
fondamental de la th\'{e}orie de l'\'{e}limination, c'est-\`{a}-dire,
l'existence des z\'{e}ros communs d'une famille donn\'{e}e de polyn\^{o}mes.

Ce probl\`{e}me a \'{e}t\'{e} abord\'{e} par plusieurs auteurs avant Mertens
comme J.J. Sylvester ou Cayley dans le cas o\`{u} le nombre des
polyn\^{o}mes co\"{\i}ncide avec celui des variables en utilisant la notion
de r\'{e}sultant, qui est ici, une forme d'inertie de degr\'{e} z\'{e}ro.

En consid\'{e}rant des polyn\^{o}mes g\'{e}n\'{e}riques (les coefficients
sont des\ ind\'{e}termi- n\'{e}es) homog\`{e}nes,  A.Hurwitz a
\'{e}tudi\'{e} dans \cite{HRWZ} l'id\'{e}al gradu\'{e} associ\'{e}
des formes d'inertie ce qui correspond \`{a} l'homologie du complexe de
Koszul d\'{e}fini par ces polyn\^{o}mes. Cette \'{e}tude a \'{e}t\'{e}
reprise par J.P.Jouanolou \cite{JONLOU1} et \cite{JONLOU2} dans le cadre de
la th\'{e}orie des sch\'{e}mas, en utilisant notamment des m\'{e}thodes
homologiques.

\noindent Dans ce travail on introduit l'id\'{e}al des formes d'inertie relatives
\`{a} des polyn\^{o}mes $f_1,...,f_r$ g\'{e}n\'{e}riques plurihomog\`{e}nes
c'est \`{a} dire homog\`{e}nes par rapport \`{a} $s$ paquets de variables $%
\underline{X}_1,...,\underline{X}_s$, dont la partie plurihomog\`{e}ne $%
{\frak A}$ de multidegr\'{e} $\left( 0,...,0\right) $ est l'id\'{e}al
r\'{e}sultant. On donne un certain nombre de caract\'{e}risations qui nous
permettent de retrouver les formules de Perron et de Perrin donn\'{e}es dans
le cas $s=1$.

En adoptant ici les m\'{e}thodes utilis\'{e}es par J.P.Jouanolou \cite
{JONLOU2}, nous d\'{e}finissons le complexe de Koszul $K^{\bullet }$ par les
polyn\^{o}mes $\left( f_i\right) $, et le complexe de \v {C}ech ${\cal C}%
^{\bullet }$ par les mon\^{o}mes $$\left( \sigma
_{i_1,...i_s}=X_{1,i_1}...X_{s,i_s}\right) _{i_1,...,i_s}.$$ Ceci nous permet
d'\'{e}tudier deux suites spectrales $^{\prime }E$ et $^{\prime \prime }E$
associ\'{e}es au bicomplexe $K^{\bullet }\bigotimes_A{\cal C}^{\bullet }$
et conduit \`{a} l'\'{e}tude de la cohomologie de Koszul et \`{a} la
cohomologie locale, cette derni\`{e}re n'est autre que la cohomologie de \v
{C}ech.

Nous donnons enfin, un r\'{e}sultat analogue au th\'{e}or\`{e}me de Hurwitz 
\cite{HRWZ} et dans un cas particulier nous retrouvons un th\'{e}or\`{e}me
de N.H.Mcoy \cite{MCOY}.

\section{Donn\'{e}es et notations}

Soient $K$ un anneau commutatif int\`{e}gre et $\underline{X}_{1},...,%
\underline{X}_{s}$ des paquets de variables avec 
$$\underline{X}_{j}=(X_{j,1},...,X_{j,n_{j}+1})$$
 pour tout$\ 1\leq j\leq s.$ On suppose de
plus que $n_{s}\geqslant n_{s-1}\geqslant ...\geqslant n_{1}\geqslant 1$.

Soit $r$ $\in {\Bbb N}$ fix\'{e} et $d_{i,1},...,d_{i,s}\,\,\,$avec $%
i=1,\ldots ,r.$ des entiers naturels non nuls. Pour tout $i=1,\ldots ,r$, on
consid\`{e}re le polyn\^{o}me g\'{e}n\'{e}rique homog\`{e}ne par rapport
\`{a} chaque paquet de variables $\underline{X}_{j\text{ }}$de degr\'{e} $%
d_{i,j}$ donn\'{e} par

\begin{equation}
f_i=\sum  U_{i,\underline{\alpha }_1,...,\underline{\alpha }_s}{}%
\underline{X}_1^{\underline{\alpha }_1}...\underline{X}_s^{\underline{\alpha 
}_s}  \label{DN1}
\end{equation}

\noindent la sommation \'etant prise pour $\alpha _1\in {\Bbb N}^{n_1+1} \mbox{tel que}  \left| 
\underline{\alpha }_1\right| =d_{i,1},  \ldots  , \underline{\alpha }%
_s\in {\Bbb N}^{n_s+1}\text{ tel}$ $\text{ que }\left| \underline{\alpha }_s\right|
=d_{i,s}$  et o\`{u} les $U_{i,\underline{\alpha }_1,...,\underline{\alpha }_s}{}$ sont
des ind\'{e}termin\'{e}es sur $K$ et $\underline{X}_j^{\underline{\alpha }%
_j}=X_{j,1}^{\alpha _{j,1}}...X_{j,n_j+1}^{\alpha _{j,n_j+1}}.$

On d\'{e}signe par $A=K\left[ U_{i,\underline{\alpha }_1,...,\underline{%
\alpha }_s}\right] \quad (1\leqslant i\leqslant r$ $\,$et $\left| \underline{%
\alpha }_j\right| =d_{i,j}$ pour tout $\ j=1,\ldots ,s),$ l'anneau des
coefficients universels.

Pour tout $j=1,\ldots ,s$, l'alg\`{e}bre de polyn\^{o}mes $%
C_{j}=A[X_{j,1},...,X_{j,n_{j}+1}]$ est not\'{e}e $A[\underline{X}_{j}].$ Si
l'on consid\`{e}re que $\deg (X_{j,l})=1$ pour $1\leqslant l\leqslant
n_{j}+1 $, cette alg\`{e}bre est naturellement ${\Bbb N-}$gradu\'{e}e.
L'alg\`{e}bre $C=C_{1}\otimes _{A}...\otimes _{A}C_{s}=A[\underline{X}%
_{1},...,\underline{X}_{s}]$ est ${\Bbb N}^{s}-$ gradu\'{e}e par : 
\begin{equation}
\begin{array}{l}
\deg (a)=(0,...,0)\in {\Bbb N}^{s},\text{ pour tout}\ a\in A, \\ 
\deg (X_{j,l})=(0,...,0,1,0,...,0)\in {\Bbb N}^{s}\text{ ici }1\text{ est
situ\'{e} sur la }j-\text{i\`{e}me position}.
\end{array}
\label{DN2}
\end{equation}
pour tout $j=1,\ldots ,s$ et $1\leqslant l\leqslant n_{j}+1$

\begin{defi}
Tout polyn\^{o}me de $C$ homog\`{e}ne par rapport \`{a} la ${\Bbb N}^{s}-$
graduation ainsi d\'{e}finie est dit {\it polyn\^{o}me plurihomog\`{e}ne}.
\end{defi}

Les polyn\^{o}mes $f_{1},...,f_{r}$ sont donc plurihomog\`{e}nes de
degr\'{e} $\deg (f_{i})=\underline{d}_{i}=(d_{i,1},...,d_{i,s})$, en
particulier l'alg\`{e}bre quotient $B=\frac{C}{(f_{1},...,f_{r})}$ de $C$
par l'id\'{e}al $(f_{1},...,f_{r})$ engendr\'{e} par $f_{1},...,f_{r}$ est $%
{\Bbb N}^{s}-$ gradu\'{e}e.

On note par ${\frak M}$ l'id\'{e}al de $C=A[\underline{X}_1,...,\underline{X}%
_s]$ engendr\'{e} par les $q$ mon\^{o}mes : 
\begin{equation}
\sigma _{i_1...i_s}=X_{1,i_1}...X_{s,i_s}\text{ o\`{u} }(i_1,...,i_s)\in 
{\prod }_ {j=1}^{s}[1,n_j+1]\,  \label{DN3}
\end{equation}

\noindent o\`{u} $q={\prod }_{j=1}^{s}(1+n_j)$ et on notera $%
\sigma _q=X_{1,n_1+1}...X_{s,n_s+1}$.

\section{Formes d'inerties}

Comme l'alg\`{e}bre $C$ est ${\Bbb N}^{s}-$gradu\'{e}e et les polyn\^{o}mes $%
f_{1},...,f_{r}$ plurihomog\`{e}nes, alors l'alg\`{e}bre quotient $B\ $est $%
{\Bbb N}^{s}-$ gradu\'{e}e. Notons par $B_{\sigma _{i_{1}...i_{s}}}$ le
localis\'{e} de $B$ par $\sigma _{i_{1}...i_{s}}$ muni de la ${\Bbb Z}^{s}-$
graduation provenant de la ${\Bbb N}^{s}-$ graduation de $B.$ Alors la
surjection canonique $p:C\longrightarrow B$ et le morphisme de $A-$%
alg\`{e}bres $\pi :b\longmapsto (\frac{b}{1},...,\frac{b}{1})$ de $B$ \`{a}
valeurs dans $\prod\limits_{i_{1},...,i_{s}}B_{\sigma _{i_{1}...i_{s}}}$
sont gradu\'{e}s de degr\'{e} $(0,...,0)\in {\Bbb N}^{s}.$ On a donc 
\[
\ker \pi =\{b\in B\mid \forall m\in {\frak M},\exists \nu \in {\Bbb N}%
,m^{\nu }b=0\}. 
\]

\begin{defi}
L'image r\'{e}ciproque ${\cal T=}p^{-1}\left( \ker \pi \right) $ est
appel\'{e}e l'{\it id\'{e}al des formes d'inertie} des polyn\^{o}mes $%
f_{1},...,f_{r}$, et la partie plurihomog\`{e}ne ${\cal T}_{(0,...,0)}={\cal %
T}$ $\cap A={\frak A}$ de degr\'{e} $(0,...,0)$, est appel\'{e}e l' {\it %
id\'{e}al r\'{e}sultant} de $f_{1},...,f_{r}$.
\end{defi}

Les formes d'inertie, c'est-\`{a}-dire les polyn\^{o}mes de ${\cal T}$ sont
caract\'{e}ris\'{e}es par la proposition suivante :

\begin{propo}
Pour tout $f\in C$ les propri\'{e}t\'{e}s suivantes sont \'{e}quivalentes :

1.$\;f\ $est une forme d'inertie.

2. Il existe $\nu \in {\Bbb N}$ tel que $\sigma _{_{i_{1}...i_{s}}}^{\nu }f$
est dans l'id\'{e}al engendr\'{e} par $f_{1},...,f_{r},$ quels que soient $%
i_{1},...,i_{s}$.

3. Il existe $\nu \in {\Bbb N}$ et il existe $i_{1},...,i_{s}$ tels que $%
\sigma _{_{i_{1}...i_{s}}}^{\nu }f$ est dans l'id\'{e}al engendr\'{e} par $%
f_{1},...,f_{r}.$

4. Il existe des entiers naturels $\nu _{1}$,...,$\nu _{s}\ $tels que, pour
tout $(\underline{\alpha }_{1},...,\underline{\alpha }_{s})\in
\prod\limits_{j=1}^{s}{\Bbb N}^{n_{j}+1}$ avec $\left| \underline{\alpha }%
_{j}\right| =$ $\nu _{j}$ $(1\leqslant j\leqslant s),$ le polyn\^{o}me $%
\underline{X}_{1}^{\underline{\alpha }_{1}}...\underline{X}_{s}^{%
\underline{\alpha }_{s}}f$ $\ $est dans l'id\'{e}al engendr\'{e} par $%
f_{1},...,f_{r}.$
\end{propo}

\medskip

On \'{e}crit $\sigma _{_{i_{1}...i_{s}}}^{\nu }f=0$ dans $B$ pour exprimer
que $\sigma _{_{i_{1}...i_{s}}}^{\nu }f$ est dans l'id\'{e}al engendr\'{e}
par $f_{1},...,f_{r}.$

La d\'{e}monstration de cette proposition repose sur le lemme suivant :

\begin{lemm}
Pour tout $i_{1},\ldots ,i_{s}\ $on a $\ker \pi =\ker (can:B\longrightarrow
B_{\sigma _{i_{1}\ldots i_{s}}}),$ o\`{u} $can$ d\'{e}signe la projection
canonique de $B$ sur $B_{\sigma _{i_{1}\ldots i_{s}}}$ .
\end{lemm}

\noindent {\it D\'{e}monstration}{\bf . }Il suffit de montrer que 
\[
\ker (can:B\longrightarrow B_{\sigma _{i_{1}\ldots i_{s}}})=\ker
(can:B\longrightarrow B_{\sigma _{j_{1}\ldots j_{s}}}) 
\]
pour tous mon\^{o}mes $\sigma _{i_{1}\ldots i_{s}}$ et $\sigma _{j_{1}\ldots
j_{s}}$ (\ref{DN3}). En effet, le diagramme

\[
\begin{array}{lll}
\,\,\,\,\,\,\,\,\,\,\,\,\,\,B & \stackrel{can}{\longrightarrow } & B_{\sigma
_{i_{1}\ldots i_{s}}} \\ 
~~can\downarrow &  & \,\downarrow \gamma _{2} \\ 
\,\,\,\,\,\,\,\,\,\,\,\,\,\,B_{\sigma _{j_{1}\ldots j_{s}}} & \stackrel{%
\gamma _{1}}{\longrightarrow } & B_{\sigma _{i_{1}\ldots i_{s}}\sigma
_{j_{1}\ldots j_{s}}}
\end{array}
\]
est commutatif, et, puisque $\sigma _{j_{1}\ldots j_{s}}$ n'est pas un
diviseur de z\'{e}ro dans $B_{\sigma _{i_{1}\ldots i_{s}}},$ alors les
homomorphismes $\gamma _{1}$ et $\gamma _{2}$ sont injectifs, par
cons\'{e}quent on a $\ker (can:B\longrightarrow B_{\sigma _{i_{1}\ldots
i_{s}}})=\ker (can:B\longrightarrow B_{\sigma _{j_{1}\ldots j_{s}}}).$
D'o\`{u}$\ $\ le lemme.%
\endproof%
%

\begin{rmk}
Il r\'{e}sulte aussit\^{o}t de la caract\'{e}risation de l'id\'{e}al des
formes d'inertie et du fait que les polyn\^{o}mes sont plurihomog\`{e}nes,
que l'id\'{e}al ${\cal T}$ est ${\Bbb N}^{s}-$gradu\'{e}.
\end{rmk}

Soit $i_1,...,i_s$ tels que $1\leq i_j\leq n_j+1.$ Pour tous $i=1,\ldots ,r$
et $j=1,\ldots ,s,$ on d\'{e}signe par $U_{i,i_1,...,i_s}$ le coefficient de 
$X_{1,i_1}^{d_{i,1}}...X_{s,i_s}^{d_{i,s}}$ dans $f_i,$ et on notera dans
toute la suite : 
\[
\begin{array}{ll}
\varepsilon _i & =U_{i,n_1+1,...,n_s+1} \\ 
\tau _i & =X_{1,n_1+1}^{d_{i,1}}...X_{s,n_s+1}^{d_{i,s}} \\ 
h_i & =f_i-\varepsilon _i\tau _i \\ 
\,\underline{\widetilde{X}}_j & =(X_{j,1},...,X_{j,n_j},1)\text{ }(\text{on
substitue }1\text{ \`{a} }X_{j,n_j+1}\text{)} \\ 
\widetilde{f} & =f(\,\,\,\underline{\widetilde{X}}_1\,,...,\underline{%
\widetilde{X}}_s),
\end{array}
\]
o\`{u} $f$ est un polyn\^{o}me de $C=A[\underline{X}_1,...,\underline{X}_s],$
que l'on regardera dans la suite sous la forme suivante $f=f(\,\varepsilon
_1,...,\varepsilon _r,\underline{X}_1,...,\underline{X}_s)\in A^{\prime
}[\,\varepsilon _1,...,\varepsilon _r,\underline{X}_1,...,\underline{X}_s]$
et $A^{\prime }=K\left[ U_{i,\underline{\alpha }_1,...,\underline{\alpha }_s}%
\right] $ avec $U_{i,\underline{\alpha }_1,...,\underline{\alpha }_s}\neq
\varepsilon _i\,(1\leqslant i\leqslant r).$ Dans ces nouvelles notations, on
a une deuxi\`{e}me caract\'{e}risation des formes d'inertie:

\begin{propo}
Pour tout polyn\^{o}me $f$ plurihomog\`{e}ne dans $C$ de degr\'{e} $(\nu
_1,...,\nu _s)$ les assertions suivantes sont \'{e}quivalentes :

\begin{enumerate}
\item  $f\in {\cal T}$

\item  $\frac{f}{X_{1,n_{1}+1}^{\nu _{1}}...X_{s,ns+1}^{\nu _{s}}}=0$ dans $%
\left( B_{\sigma _{q}}\right) _{(0,\ldots ,0)}$,

\item  $f(-\widetilde{h}_1,...,-\widetilde{h}_r,\underline{\widetilde{X}}%
_1,...,\widetilde{\underline{X}}_s)=0$ dans $C.$
\end{enumerate}
\end{propo}

Montrons tout d'abord le lemme suivant :

\begin{lemm}
Pour tout $i_1\ldots i_s,$ il existe un isomorphisme de $A^{\prime }[%
\underline{X}_1,...,\underline{X}_s]-$alg\`{e}bre: 
\[
B_{\sigma _{i_1\ldots i_s}}\longrightarrow A^{\prime }[\underline{X}_1,...,%
\underline{X}_s]_{\sigma _{i_1\ldots i_s}}. 
\]
\end{lemm}

\noindent {\it D\'{e}monstration}{\bf . }On peut supposer que $%
i_1=n_1+1,\ldots ,i_s=n_s+1,$ c'est \`{a} dire $\sigma _{i_1\ldots
i_s}=\sigma _q.$

Soit donc $\varphi $ l'homomorphisme de $A^{\prime }[\underline{X}_1,...,%
\underline{X}_s]-$alg\`{e}bres

\[
\varphi :C=A^{\prime }[\underline{X}_1,...,\underline{X}_s][\varepsilon _1
,...,\varepsilon _r]\longrightarrow A^{\prime }[\underline{X}_1,...,%
\underline{X}_s]_{\sigma _q} 
\]
d\'{e}fini par $\varphi \left( \,\varepsilon _i\right) =-\frac{h_i}{\tau _i}%
, $ o\`{u} $\,\,\,\tau _i=X_{1,n_1+1}^{d_{i,1}}...X_{s,n_s+1}^{d_{i,s}}\ $et 
$h_i=f_i-\varepsilon _i\tau _i.$

L'homomorphisme $\varphi \,$est bien d\'{e}fini, car$\,\,\tau _i$ est
inversible dans $A^{\prime }[\underline{X}_1,...,\underline{X}_s]_{\sigma
_q},$ et, puisque $\varphi (f_i)=0,$ pour tout $i=1,...,r,$ et $\varphi
(\sigma _q)$ est inversible dans $A^{\prime }[\underline{X}_1,...,\underline{%
X}_s]_{\sigma _q},$ alors $\varphi $ induit un homomorphisme de $A^{\prime }[%
\underline{X}_1,...,\underline{X}_s]-$alg\`{e}bres $\psi :B_{\sigma
_q}\longrightarrow A^{\prime }[\underline{X}_1,...,\underline{X}_s]_{\sigma
_q}.$ Soit $\psi ^{\prime }$ l'homomorphisme compos\'{e} 
\[
\psi ^{\prime }:A^{\prime }[\underline{X}_1,...,\underline{X}_s]_{\sigma _q}%
\stackrel{j}{\longrightarrow }C_{\sigma _q}\stackrel{p}{\longrightarrow }%
B_{\sigma _q} 
\]

o\`{u} $j$ est l'injection canonique et $p$ et l'homomorphisme induit de la
surjection canonique de $C$ dans $B=\frac C{(f_1,...,f_r)}.$ On d\'{e}duit
que $\psi \circ \psi ^{\prime }=Id_{A^{\prime }[\underline{X}_1,...,%
\underline{X}_s]_{\sigma _q}},$ de m\^{e}me on a $\psi ^{\prime }\circ \psi
=Id_{B_{\sigma _q}}$.%
\endproof%
%

\noindent {\it D\'{e}monstration de la proposition}{\bf .} Il r\'{e}sulte de
la premi\`{e}re caract\'{e}risation des formes d'inertie que les deux
premi\`{e}res assertions sont \'{e}quivalentes.

Soit $f\in {\cal T}$, plurihomog\`{e}ne. de degr\'{e} $(\nu _1,...,\nu _s).$
On a 
$$f=f(\,\varepsilon _1,...,\varepsilon _r,\underline{X}_1,...,\underline{%
X}_s)=0$$
dans $B_{\sigma _q}$ (premi\`{e}re caract\'{e}risation des formes
d'inertie), donc 
\[
\psi (f)=f(-\frac{h_1}{\tau _1},...,-\frac{h_r}{\tau _r},\underline{X}_1,...,%
\underline{X}_s)=0\text{ dans }A^{\prime }[\underline{X}_1,...,\underline{X}%
_s]_{\sigma _q}, 
\]
et, si on remplace $X_{j,n_j+1}$ par $1$, pour tout $j=1,...,s,$ dans $\psi
(f),$ on en d\'{e}duit que 
$$f(-\widetilde{h}_1,...,-\widetilde{h}_r,%
\widetilde{\underline{X}}_1,...,\underline{\widetilde{X}}_s)=0$$
 dans $C,$
d'o\`{u} 1$\Longrightarrow 3.$

R\'{e}ciproque. Supposons que $f(-\widetilde{h}_{1},...,-\widetilde{h}_{r},%
\widetilde{\underline{X}}_{1},...,\underline{\widetilde{X}}_{s})=0$ dans $C$%
. Par homog\'{e}n\'{e}isation par rapport \`{a} chaque paquet de variables $%
\underline{X}_{j}$ \`{a} l'aide de la variable $X_{j,n_{j}+1}$, on obtient

\[
f(-\frac{h_{1}}{\tau _{1}},...,-\frac{h_{r}}{\tau _{r}},\frac{\underline{X}%
_{1}}{X_{1,n1+1}},...,\frac{\underline{X}_{s}}{X_{s,ns+1}})=0\text{ dans }%
C_{\sigma _{q}} 
\]
car $h_{i}$ est plurihomog\`{e}ne. de degr\'{e} $\underline{d}%
_{i}=(d_{i,1},...,d_{i,s})$. Or, pour tout $i=1,\ldots ,r$, $-\frac{h_{i}}{%
\tau _{i}}=\varepsilon _{i}$ dans $B_{\sigma _{q}},$ et comme $f$ est
plurihomog\`{e}ne, on en d\'{e}duit que $\frac{f}{1}=0$ dans $B_{\sigma
_{q}},$ d'o\`{u} $3\Longrightarrow 1.$%
\endproof%
%

\medskip

On va donner ici une autre caract\'{e}risation de l'id\'{e}al r\'{e}sultant $%
{\frak A}$ des polyn\^{o}mes $f_1,...,f_r$ c'est \`{a} dire, des formes
d'inertie de degr\'{e} $\left( 0,\ldots ,0\right) $.

Consid\'{e}rons l'homomorphisme de $C-$alg\`{e}bres $F:T_i\longmapsto \frac{%
f_i}{\tau _i}$ de $C[T_1,...,T_r]$ \`{a} valeurs dans$\ C_{\sigma _q}$ et
l'automorphisme $G:\,\varepsilon _i\longmapsto \varepsilon _i-T_i$ de
l'alg\`{e}bre $C[T_1,...,T_r]$ sur l'anneau $A^{\prime }[\underline{X}_1,...,%
\underline{X}_s,T_1,...,T_r],$o\`{u} $T_1,...,T_r$ sont des nouvelles
variables.

\begin{propo}
$KerF$ co\"{\i}ncide avec l'image par $G$ de l'id\'{e}al ${\cal T}%
[T_{1},...,T_{r}]$ .
\end{propo}

\noindent {\it D\'{e}monstration}{\bf . }On consid\`ere l'homomorphisme de $A^{\prime
}[\underline{X}_1,...,\underline{X}_s]_{\sigma _q}-$ alg\`{e}bres de  $%
C_{\sigma _q}[T_1,...,T_r]$ \`{a} valeurs dans $C_{\sigma _q}$ d\'{e}fini
par : 
\[
T_i\longmapsto \frac{f_i}{\tau _i}\text{ et }\,\varepsilon _i\longmapsto -%
\frac{h_i}{\tau _i}\text{ .} 
\]
Puisque, pour tout $1\leqslant i\leqslant r,$ on a $-\frac{h_i}{\tau _i}%
=\,\varepsilon _i$ dans $B_{\sigma _q},$ alors l'id\'{e}al engendr\'{e} par $%
f_1,...,f_r$ est le noyau de l'homomorphisme d\'{e}fini ci dessus, qui
induit par cons\'{e}quent, un isomorphisme $\omega :B_{\sigma
_q}[T_1,...,T_r]\longrightarrow C_{\sigma _q}.$ On sait que l'homomorphisme
canonique d'anneaux $u:C\longrightarrow B_{\sigma _q}$ admet pour noyau
l'id\'{e}al ${\cal T}$ des formes d'inertie, et ${\cal T}[T_1,...,T_r]$ est
le noyau du morphisme d'alg\`{e}bres $\overline{u},$ qui prolonge $u$ \`{a} $%
C[T_1...,T_r].\ $Le diagramme commutatif suivant

\begin{center}
$
\begin{array}{lll}
C[T_1,...,T_r] & \stackrel{\overline{u}}{\longrightarrow } & B_{\sigma
_q}[T_1,...,T_r] \\ 
\,\,\,\,\,\,\,\,\,\,\,\,G\downarrow &  & \downarrow \omega \\ 
C[T_1,...,T_r] & \stackrel{F}{\longrightarrow } & C_{\sigma _q}
\end{array}
$
\end{center}

\noindent montre que $\ker F=G({\cal T}[T_1,...,T_r])$. 
\endproof%
%

\noindent On d\'{e}duit :

\begin{coro}
{\bf Formule de Perron. }L'homomorphisme de $A$- alg\`{e}bres $\widetilde{F}%
:T_{i}\longmapsto \widetilde{f_{i}},$ de $A[T_{1},...,T_{r}]$ \`{a} valeurs
dans $A[{\widetilde{\underline{X}_{1}}},...,{%
\widetilde{\underline{X}_{s}}}]$ a pour noyau : $\ker \widetilde{F}=G({\cal T})[T_{1},...,T_{r}].$
\end{coro}

\begin{coro}
{\bf Formule de Perrin. }Pour tout $a\in A$, les propri\'{e}t\'{e}s
suivantes sont \'{e}quivalentes

$\ i)$ $a\in {\frak A,}$

ii) Il existe un polyn\^{o}me $Q\in A[T_{1},...,T_{r}]$ sans terme constant
tel que $a=Q(\widetilde{f}_{1},...,\widetilde{f}_{r})$
\end{coro}

\noindent {\it D\'{e}monstration}{\bf . }Soit $a=a(\varepsilon
_1,...,\varepsilon _r)\in {\frak A}$ consid\'{e}r\'{e} comme polyn\^{o}me en 
$\varepsilon _1,...,\varepsilon _r$. L'\'{e}galit\'{e} $\ker \widetilde{F}=G(%
{\frak A}[T_1,...,T_r])$ montre que $G(a)=$ $a(\varepsilon
_1-T_1,...,\varepsilon _r-T_r)\in \ker \widetilde{F}$.

En appliquant $\widetilde{F}$ au polyn\^{o}me $Q(T_1,...,T_r)=a(\varepsilon
_1,...,\varepsilon _r)-a(\varepsilon _1-T_1,...,\varepsilon _r-T_r)$, on
obtient $a(\varepsilon _1,...,\varepsilon _r)=Q(\widetilde{f}_1,...,%
\widetilde{f}_r)$.

R\'{e}ciproquement, s'il existe $Q(T_1,...,T_r)\in $ $A[T_1,...,T_r]$ sans
terme constant tel que $a=Q(\widetilde{f}_1,...,\widetilde{f}_r)$, alors $%
a-Q(T_1,...,T_r)\in \ker \widetilde{F}=G({\frak A}[T_1,...,T_r])$, donc on
peut \'{e}crire $a-Q(T_1,...,T_r)=\sum\limits_\alpha a_\alpha (\varepsilon
_1-T_1)^{\alpha _1}...(\varepsilon _r-T_r)^{\alpha _r}$ avec $a_\alpha \in 
{\frak A}$. D'o\`{u}, par sp\'{e}cialisation $T_i\longmapsto 0$ pour tout $%
i=1,\ldots ,r,$ on d\'{e}duit $a\in {\frak A}$. 
\endproof%
%

\section{Complexes de Koszul et de \v {C}ech}

On \'{e}tudie ici le complexe de Koszul et la cohomologie locale
associ\'{e}s \`{a} des polyn\^{o}mes plurihomog\`{e}nes.

Dans ce paragraphe, on d\'{e}signe par $A$ un anneau commutatif, $M$ un $A-$
module libre et $\underline{a}=(a_1,...,a_r)$ une suite d'\'{e}l\'{e}ments
de $A$.

\subsection{Complexe de Koszul}

Soient $(e_1,...,e_r)$ la base canonique du $A-$ module libre $A^r$. Pour $%
0\leq p\leq r$, la $p-i\grave{e}me$ puissance ext\'{e}rieure $%
\bigwedge^p(A^r)$ du $A-$module $A^r$ est un $A-$module libre de rang $C_r^p=%
\frac{r!}{p!\left( r-p\right) !}$ et de base $\left( e_{i_1}\wedge ...\wedge
e_{i_p}\right) _{1\leq i_1<...<i_p\leq r},$ avec la convention $%
\bigwedge^0(A^r)=A^r$.

Le complexe de cha\^{i}nes de Koszul associ\'{e} \`{a} la suite $\underline{a%
},$ not\'{e} $(K_{\bullet }(\underline{a},A),d_{\bullet }),$ est d\'{e}fini
par 
\[
K_{p}(\underline{a},A)=\bigwedge\nolimits^{p}(A^{r}) 
\]
pour tout $p=0,\ldots ,r,$ \thinspace telle que pour $1\leq
i_{1}<...<i_{p}\leq r,$ et 
\[
d_{p}:K_{p}(\underline{a},A)\longrightarrow K_{p-1}(\underline{a},A) 
\]
est donn\'{e}e par 
\[
d_{p}(e_{i_{1}}\wedge ...\wedge
e_{i_{p}})=\sum_{k=1}^{p}(-1)^{k+1}a_{i_{k}}e_{i_{1}}\wedge ...\wedge 
\widehat{e_{i_{k}}}\wedge ...\wedge e_{i_{p}}. 
\]
On d\'{e}finit aussi, pour tout $A-$module $M$, le complexe de cha\^{\i}nes
de Koszul associ\'{e} \`{a} la suite $\underline{a}$ et au module $M$ par : 
\begin{equation}
K_{\bullet }(\underline{a},M)=K_{\bullet }(\underline{a},A)\bigotimes%
\nolimits_{A}M.  \label{Comp Kosz1}
\end{equation}
Si $A$ est ${\Bbb N-}$gradu\'{e} alors le $A-$module $\bigwedge^{p}(A^{r})$
est ${\Bbb N-}$gradu\'{e} par : 
\begin{equation}
\left\{ 
\begin{array}{l}
\deg (e_{k})=\nu _{k},\text{ pour }1\leq k\leq r \\ 
\deg (e_{i_{1}}\wedge ...\wedge e_{i_{p}})=\nu _{i_{1}}+...+\nu _{i_{p}},%
\text{ pour }1\leq i_{1}<...<i_{p}\leq r\text{ }
\end{array}
\right.  \label{Kosz-grad}
\end{equation}

Si $a_{1},...,a_{r}$ sont homog\`{e}nes de degr\'{e}s $\nu _{1},...,\nu _{r}$
alors les diff\'{e}rentielles sont homog\`{e}nes de degr\'{e} $0$. On
obtient ainsi une graduation du complexe de Koszul $K_{\bullet }(\underline{a%
},M)$, et pour tout $p=1,\ldots ,r$, on a un isomorphisme $K_{p}(\underline{a%
},A)\longrightarrow \bigoplus\limits_{1\leq i_{1}<\ldots <i_{p}\leq r}A\left[
-\nu _{i_{1}}-\ldots -\nu _{i_{p}}\right] ,$ o\`{u} $A\left[ \nu \right] $
est le $A-$module gradu\'{e} obtenu par d\'{e}calage de la graduation de $%
\nu \in {\Bbb Z}$, c'est \`{a} dire, $(A\left[ \nu \right] )_{t}=A_{\nu +t}$.

\subsection{Complexe de \v {C}ech}

Pour tout $I=\left\{ i_1,\ldots ,i_p\mid i_1<...<i_p\right\} \subset
\{1,\ldots ,r\},$ on d\'{e}signe par $M_{a_I}$ le $A-$module $%
M_{a_{i_1}\ldots a_{i_p}}$ localis\'{e} de $M$ par le produit $%
a_{i_1}...a_{i_p}$.

Le complexe de \v{C}ech associ\'{e} \`{a} la suite $\underline{a}=\left(
a_{1},...,a_{r}\right) $ et au module $M$ est le complexe not\'{e} $\left( 
\stackrel{\vee }{C^{\bullet }}(\underline{a},M),d^{\bullet }\right) $
d\'{e}fini par

\begin{equation}
\stackrel{\vee }{C^{p}}(\underline{a},M)=\prod\limits_{\mid I\mid
=p+1}M_{a_{I}}  \label{Comp.Cech}
\end{equation}
et les diff\'{e}rentielles $d^{p}:\,\stackrel{\vee }{C^{p}}(\underline{a}%
,M)\longrightarrow \stackrel{\vee }{C^{p+1}}(\underline{a},M)$ o\`{u} pour
tout $0\leq p\leq r-1$ sont d\'{e}finies, pour tout $m=(m_{I})_{\mid I\mid
=p+1}$ et pour tout $J\subset \{1,...,r\}$ de cardinal $p+2,$ par :

\begin{equation}
\left( d^p(m)\right) _J=\sum\limits_{j\in J}(-1)^{l_j-1}\frac{m_{J-\{j\}}}1
\label{Dif. Cech}
\end{equation}
o\`{u} $l_j$ est la position de $j$ dans $J=\left\{ j_1,\ldots ,j_{p+2}\mid
j_1<\ldots <j_{p+2}\right\} $ et $\frac{m_{J-\{j\}}}1$ est l'image de $%
m_{J-\{j\}}$ par le morphisme canonique $M_{a_{J-\{j\}}}\longrightarrow
(M_{a_{J-\{j\}}})_{a_j}=M_{a_J\ }.$

On note par $\stackrel{\vee }{H^{\bullet }}(\underline{a},M)$ la cohomologie
du complexe du \v {C}ech$\ \stackrel{\vee }{C^{\bullet }}(\underline{a},M).$

Soit maintenant $\alpha :m\longmapsto \left( \frac m1,...,\frac m1\right) $
l'homomorphisme de $A-$modules de $M$ \`{a} valeurs dans $%
\prod\limits_{i=1}^rM_{a_i}.$ On v\'{e}rifie ais\'{e}ment que $d^0\circ
\alpha =0$ (voir \ref{Dif. Cech}). On a donc un autre complexe : 
\begin{equation}
0\rightarrow M\stackrel{\alpha }{\rightarrow }\prod\limits_{i=1}^rM_{a_i}%
\stackrel{d^0}{\rightarrow }\prod\limits_{\mid I\mid =2}M_{a_I}\stackrel{d^1%
}{\rightarrow }\ldots \stackrel{d^{r-2}}{\rightarrow }M_{a_1...a_r}\stackrel{%
d^{r-1}}{\rightarrow }0,  \label{CechAug1}
\end{equation}
appel\'{e} complexe de \v {C}ech augment\'{e}, qu'on d\'{e}signera par $%
\stackrel{\vee }{{\cal C}^{\bullet }}\left( \underline{a},M\right) ;$
d'o\`{u} les relations :

\[
H^0\left( \stackrel{\vee }{{\cal C}^{\bullet }}(\underline{a},M)\right)
=\ker \alpha ,\text{ }H^1\left( \stackrel{\vee }{{\cal C}^{\bullet }}(%
\underline{a},M)\right) =\frac{\stackrel{\vee }{H^0}(\underline{a},M)}{\Im}\alpha ,
\]
et$,$ 
\[
H^i\left( \stackrel{\vee }{{\cal C}^{\bullet }}(\underline{a},M)\right) =%
\stackrel{\vee }{H^{i-1}}(\underline{a},M)\quad \text{pour }i\geq 2. 
\]

Soit $J=\left( a_1,...,a_r\right) $ l'id\'{e}al de $A$ engendr\'{e} par $%
a_1,...,a_r$. On note 
$$H^{\bullet }\left( \stackrel{\vee }{{\cal C}^{\bullet
}}(\underline{a},M)\right) =H_J^{\bullet }\left( M\right) .$$

Consid\'{e}rons le sch\'{e}ma affine $X=spec(A)$ et $\,U=\bigcup%
\limits_{i=1}^rD(a_i)$ la r\'{e}union des ouverts affines $D(a_i)$ de $X,$
d\'{e}finis respectivement par $a_1,...,a_r.$

Le complexe de \v{C}ech $\stackrel{\vee }{C}^{\bullet }\left( \underline{a}%
,M\right) $ associ\'{e} \`{a} la suite $\underline{a}=\left(
a_{1},...,a_{r}\right) $ et au module $M$, n'est rien d'autre que le
complexe de \v{C}ech habituel associ\'{e} au recouvrement $U=\left(
D(a_{i})\right) _{1\leqslant i\leqslant r}$ et\thinspace au faisceau $%
\widetilde{M}$ associ\'{e} au $A-$module $M$ (voir \cite{HRTSHORN},
Chapitre. III). D'apr\`{e}s le th\'{e}or\`{e}me de Cartan-Leray on a $%
H^{i}(U,\widetilde{M})=\stackrel{\vee }{H^{i}}(\underline{a},M),$ pour tout $%
i=0,\ldots ,r$

\begin{propo}
Soit $Y$ le sous sch\'{e}ma ferm\'{e} du sch\'{e}ma affine $X=spec\left(
A\right) $ d\'{e}fini par l'id\'{e}al $J=\left( a_1,...,a_r\right) $ de $A$.
On a alors un isomorphisme 
$$H_J^{\bullet }(M)\longrightarrow H_Y^{\bullet
}(X,\widetilde{M})$$
o\`{u} $H_Y^{\bullet }(X,\widetilde{M})$ est la cohomologie \`{a} support
dans $Y$.
\end{propo}

\noindent {\it D\'{e}monstration}. De la suite exacte longue de cohomologie 
\begin{eqnarray*}
0 &\rightarrow &H_{Y}^{0}(X,\widetilde{M})\rightarrow H^{0}(X,\widetilde{M}%
)\rightarrow H^{0}(U,\widetilde{M})\rightarrow H_{Y}^{1}(X,\widetilde{M}%
)\rightarrow \ldots \\
\ldots &\rightarrow &H_{Y}^{i}(X,\widetilde{M})\rightarrow H^{i}(X,%
\widetilde{M})\rightarrow H^{i}(U,M)\rightarrow H_{Y}^{i+1}(X,\widetilde{M}%
)\rightarrow \ldots
\end{eqnarray*}
et\thinspace du\thinspace th\'{e}or\`{e}me\thinspace de\thinspace
Cartan-Serre$\,(H^{i}(X,\widetilde{M})=0$\thinspace $\,$pour\thinspace $%
i\geqslant 1\,$\cite{GRTDK-DDN 3}) on d\'{e}duit qu'on a un isomorphisme $%
H^{i}(U,M)\simeq H_{Y}^{i+1}(X,\widetilde{M}),$ pour $i\geqslant 1.$ La
suite 
\[
0\rightarrow H_{Y}^{0}(X,\widetilde{M})\rightarrow H^{0}(X,\widetilde{M}%
)\rightarrow H^{0}(U,\widetilde{M})\rightarrow H_{Y}^{1}(X,\widetilde{M}%
)\rightarrow 0. 
\]
est exacte, $U$ \'{e}tant $X-Y.$

Il r\'{e}sulte du th\'{e}or\`{e}me de Cartan-Leray et de la d\'{e}finition
du complexe de \v{C}ech augment\'{e} que, pour $i\geqslant 2$, on a un
isomorphisme $H_{J}^{i}(M)\longrightarrow H_{Y}^{i}(X,\widetilde{M}).$ Pour $%
i=0$, on a $H_{Y}^{0}(X,\widetilde{M})=\ker (M\rightarrow \Gamma (U,%
\widetilde{M}))=\ker (\alpha )=H_{J}^{0}(M).$ On d\'{e}duit de la suite
exacte pr\'{e}c\'{e}dente que pour $i=1$ on a $H_{Y}^{1}(X,\widetilde{M}%
)=H_{J}^{1}(M),$ ce qui ach\`{e}ve la d\'{e}monstration. 
\endproof%
%

\begin{exemple}
Soient $R=K[X_1,...,X_n]$ l'anneau des polyn\^{o}mes \`{a} coefficients dans
un corps commutatif $K,$ $U$ l'ouvert $\bigcup\limits_{i=1}^nD\left(
X_i\right) $ de $spec\left( R\right) $ et soit $J=\left( X_1,...,X_n\right) $
l'id\'{e}al de $R$ engendr\'{e} par $X_1,...,X_n$. Puisque $X_1,...,X_n$ est
une $R-$suite, alors $prof(J)\geqslant n$, et par cons\'{e}quent les groupes
de cohomologie $H_J^i\left( R\right) $ sont nuls sauf pour $i\neq n$. Et
d'apr\`{e}s \ref{CechAug1}, on a $H_J^n\left( R\right) $ est le $n-i\grave{e}%
me\ $groupe de cohomologie du complexe $\stackrel{\vee }{{\cal C}^{\bullet }}%
\left( \underline{X},R\right) ,$ d'o\`{u} : 
\begin{equation}
H_J^n\left( R\right) =\frac{R_{X_1...X_n}}{\Im\left(
\bigoplus_{i=1}^nR_{X_1...\widehat{X_i}...X_n}\rightarrow
R_{X_1...X_n}\right) }=\frac 1{X_1...X_n}K[X_1^{-1},...,X_n^{-1}].
\label{EXPL}
\end{equation}
De plus, si $R$ est gradu\'{e} au moyen des $X_i$ $\left( \deg X_i=1\right)
, $ alors les groupes de cohomologie sont gradu\'{e}s, et on a $H_J^n\left(
R\right) _\nu =0$ si $\nu >-n$ et $H_J^n\left( R\right) _{-n}=K.$
\end{exemple}

\subsection{Suites spectrales associ\'{e}es au bicomplexe $K^{\bullet
\bullet }(\protect\underline{f},C)$.}

Soit $\underline{f}=\left( f_1,...,f_r\right) $ la suite de polyn\^{o}mes
g\'{e}n\'{e}riques plurihomog\`{e}nes (\ref{DN1}). Notons par : 
\begin{equation}
\left\{ 
\begin{array}{l}
K^0=C \\ 
K^{-l}=K_l\left( \underline{f},C\right) =\bigoplus\limits_{1\leqslant
i_1<...<i_l\leqslant r}C[-\underline{d}_{i_1}-...-\underline{d}_{i_l}],\quad 
\text{pour}\ l=0,...,r\ 
\end{array}
\right.  \label{Comp Kosz2}
\end{equation}
le complexe de Koszul d\'{e}fini sur $C=A[\underline{X}_1,...,\underline{X}%
_s]$ par la suite $\underline{f}$ $.$

Les polyn\^{o}mes $f_i$ sont homog\`{e}nes par la ${\Bbb N}^s-$ graduation
d\'{e}finie dans \ref{DN2}, ce qui induit une ${\Bbb N}^s-$ graduation sur $%
K^{\bullet }$.

Comme les diff\'{e}rentielles $d^{-l}:K^{-l}\longrightarrow K^{-l+1}$ sont
homog\`{e}nes de degr\'{e} $\left( 0,...,0\right) \in {\Bbb N}^s,$ on
d\'{e}duit que les groupes de cohomologie $(H^i(K^{\bullet }))_{-r\leqslant
i\leqslant 0}$ sont aussi ${\Bbb N}^s-$ gradu\'{e}s.

On d\'{e}finit maintenant le complexe de \v {C}ech augment\'{e} sur l'anneau 
$C$ par la suite de mon\^{o}mes $\underline{\sigma }=\left( \sigma
_{i_1\ldots i_s}\right) _{i_1\ldots i_s}$ (\ref{CechAug1} et \ref{DN3}) que
l'on note

\begin{equation}
\stackrel{\vee }{{\cal C}^p}=\stackrel{\vee }{{\cal C}^p}\left( \underline{%
\sigma },C\right) ,  \label{CechAug2}
\end{equation}
pour $p=1,...,q.$ On a donc un bicomplexe 
\[
K^{\bullet \bullet }=K^{\bullet \bullet }(\underline{f},C)=K^{\bullet
}\bigotimes\nolimits_C\stackrel{\vee }{{\cal C}^{\bullet }} 
\]
qu'on se propose d'\'{e}tudier dans le paragraphe suivant.

\subsubsection{Propri\'{e}t\'{e}s du bicomplexe $K^{\bullet \bullet }(%
\protect\underline{f},C)$}

Pour tous $l={\Bbb -}r,...,0$ et $p=1,...,q,$ on pose $K^{l,p}=K^l\bigotimes%
\nolimits_C\stackrel{\vee }{{\cal C}^p.\text{ }}$Consid\'{e}rant $p$ comme
indice de ligne et $l$ comme indice de colonne, on obtient le diagramme
commutatif suivant

\begin{equation}
\begin{array}{lllllllllllll}
&  & 0 &  & 0 &  &  &  & 0 &  & 0 &  &  \\ 
&  & \downarrow &  & \downarrow &  &  &  & \downarrow &  & \downarrow &  & 
\\ 
0 & \rightarrow & K^{-r,0} & \stackrel{d^{^{\prime }}}{\rightarrow } & 
K^{-r+1,0} & \stackrel{d^{^{\prime }}}{\rightarrow } & ... & \stackrel{%
d^{^{\prime }}}{\rightarrow } & K^{-1,0} & \stackrel{d^{^{\prime }}}{%
\rightarrow } & K^{0,0} & \rightarrow & 0 \\ 
&  & \downarrow d^{^{\prime \prime }} &  & \downarrow d^{^{\prime \prime }}
&  &  &  & \downarrow d^{^{\prime \prime }} &  & \downarrow d^{^{\prime
\prime }} &  &  \\ 
0 & \rightarrow & K^{-r,1} & \stackrel{d^{^{\prime }}}{\rightarrow } & 
K^{-r+1,1} & \stackrel{d^{^{\prime }}}{\rightarrow } &  & \stackrel{%
d^{^{\prime }}}{\rightarrow } & K^{-1,1} & \stackrel{d^{^{\prime }}}{%
\rightarrow } & K^{0,1} & \rightarrow & 0 \\ 
&  & \downarrow d^{^{\prime \prime }} &  & \downarrow d^{^{\prime \prime }}
&  &  &  & \downarrow d^{^{\prime \prime }} &  & \downarrow d^{^{\prime
\prime }} &  &  \\ 
&  & \vdots &  & \vdots &  &  &  & \vdots &  & \vdots &  &  \\ 
&  & \downarrow d^{^{\prime \prime }} &  & \downarrow d^{^{\prime \prime }}
&  &  &  & \downarrow d^{^{\prime \prime }} &  & \downarrow d^{^{\prime
\prime }} &  &  \\ 
0 & \rightarrow & K^{-r,q} & \stackrel{d^{^{\prime }}}{\rightarrow } & 
K^{-r+1,q} & \stackrel{d^{^{\prime }}}{\rightarrow } &  & \stackrel{%
d^{^{\prime }}}{\rightarrow } & K^{-1,q} & \stackrel{d^{^{\prime }}}{%
\rightarrow } & K^{0,q} & \rightarrow & 0 \\ 
&  & \downarrow &  & \downarrow &  &  &  & \downarrow &  & \downarrow &  & 
\\ 
&  & 0 &  & 0 &  &  &  & 0 &  & 0 &  & 
\end{array}
\label{Diag de comp}
\end{equation}
o\`{u} $d^{^{\prime }}$ (resp. $d^{\prime \prime })$ d\'{e}signe les
diff\'{e}rentielles des lignes (resp. des colonnes), obtenues par
tensorisation \`{a} partir de celles de $K^{\bullet }$ et de $\stackrel{\vee 
}{{\cal C}^{\bullet }\text{.}}$ Les complexes lignes $K^{\bullet
,p}=K^{\bullet }\bigotimes_{C}\stackrel{\vee }{{\cal C}^{p}}=K^{\bullet }(%
\underline{f},\stackrel{\vee }{{\cal C}^{p}}),$ pour $0\leqslant p\leqslant q
$ sont des complexes de Koszul${\Bbb \ N}^{s}-$ gradu\'{e}s associ\'{e}s
\`{a} la suite $\underline{f}$ et aux $C-$ modules${\Bbb \ }\stackrel{\vee }{%
{\cal C}^{p}}$. Les complexes colonnes $K^{l,\bullet }=K^{l}\bigotimes_{C}%
\stackrel{\vee }{{\cal C}^{\bullet }}=\stackrel{\vee }{{\cal C}^{\bullet }}(%
\underline{{\cal \sigma }},K^{l}),$ pour $-r\leqslant l\leqslant 0$ sont des
complexes de \v{C}ech augment\'{e}s associ\'{e}s \`{a} la suite $\underline{%
\sigma }=\left( \sigma _{i_{1}\ldots i_{s}}\right) _{i_{1}\ldots i_{s}}$ et
aux $C-$ modules $K^{l}.$

Le bicomplexe $K^{\bullet \bullet }$ donne lieu \`{a} deux suites spectrales
ayant m\^{e}me aboutissement : 
\begin{equation}
\left\{ 
\begin{array}{lllll}
^{^{\prime }}E_{1}^{l,p} & = & H^{p}(K^{l,\bullet }) & \Rightarrow & E^{l+p}
\\ 
^{^{\prime \prime }}E_{2}^{l,p} & = & H^{l}(L^{\bullet ,p}) & \Rightarrow & 
E^{l+p}
\end{array}
\right. ,\text{pour }(l,p)\in {\Bbb Z}^{2}  \label{Suite spec1}
\end{equation}
(voir\cite{GRTDK-DDN 3} \S 11, \cite{CRT-ELBRG}$,$\cite{GDMT}), o\`{u} $%
L^{p,\bullet }$ est le complexe 
\[
0\rightarrow H^{p}(K^{\bullet ,0})\rightarrow H^{p}(K^{\bullet
,1})\rightarrow ...\rightarrow H^{p}(K^{\bullet ,q})\rightarrow 0. 
\]

Puisque $K^{\bullet ,i}=K^{\bullet }\bigotimes_{C}\stackrel{\vee }{{\cal C}%
^{i}}$ est un $C-$ module plat, alors le complexe $L^{p,\bullet }$ n'est
autre que le complexe de \v{C}ech augment\'{e} associ\'{e} \`{a} la suite $%
\underline{\sigma }$ et au $C-$ module $H^{p}(K^{\bullet })$. On en
d\'{e}duit 
\begin{equation}
\left\{ 
\begin{array}{lllll}
^{^{\prime }}E_{1}^{l,p} & = & H_{{\frak M}}^{p}(K^{l}) & \Rightarrow & 
E^{l+p} \\ 
^{^{\prime \prime }}E_{2}^{l,p} & = & H_{{\frak M}}^{l}(H^{p}(K^{\bullet }))
& \Rightarrow & E^{l+p},
\end{array}
\right.  \label{Suite spec2}
\end{equation}
et, il r\'{e}sulte des d\'{e}finitions \ref{Comp Kosz2} et \ref{CechAug2},
que l'on a

\begin{center}
$\left\{ 
\begin{array}{llll}
^{^{\prime }}E_{1}^{l,p} & = & 0 & \text{pour }(l,p)\notin
\{-r,...,0\}\times \{0,...,q\} \\ 
^{^{\prime \prime }}E_{2}^{l,p} & = & 0 & \text{pour }(l,p)\notin
\{0,...,q\}\times \{-r,...,0\}.
\end{array}
\right. $
\end{center}

\subsubsection{ Supports des deux suites spectrales $^{^{\prime }}E$ et $%
^{^{\prime \prime }}E$}

Le support d'une suite spectrale $\left( E_t^{l,p}\right) \left( t\geqslant
\alpha =1\text{ ou }2\right) $ est l'ensemble des couples entiers $(l,p)$
tels que $E_\alpha ^{l,p}\neq 0$. Pour tout $j=1,...,s$, on note $U_{j\text{ 
}}=spec(C_j)-V({\frak M}_j)$ l'ouvert compl\'{e}mentaire du ferm\'{e} $V(%
{\frak M}_j)$ dans $spec(C_j)$ et $U=U_1\times _S...\times _SU_s$ l'ouvert
produit fibr\'{e} sur $S=spec\left( A\right) $ des ouverts $U_j$. Or, pour
tout $1\leqslant j\leqslant s,$ on a $U_j=\left( E_A^{n_j+1}\right) ^{*},$
donc 
\[
U=\left( E_A^{n_1+1}\right) ^{*}\times _S...\times _S\left(
E_A^{n_s+1}\right) ^{*}=spec\left( C\right) -V({\frak M}). 
\]
On d\'{e}duit la proposition suivante :

\begin{propo}
\label{première pro suite spec}On a $^{^{\prime }}E_{1}^{l,p}=0,$ pour tous $%
l\in \left\{ -r,...,0\right\} $ et $p\neq 1+\sum\limits_{j\in J}n_{j}$ pour
toute partie $J$ non vide de $\{1,...,s\}$.
\end{propo}

{\it \noindent D\'{e}monstration. }Le $C-$module $K^{l}$ est libre, donc
plat, par cons\'{e}quent on a $^{^{\prime }}E_{1}^{l,p}=H_{{\frak M}%
}^{p}(K^{l})=K^{l}\bigotimes\nolimits_{C}H_{{\frak M}}^{p}(C).$ Il suffit
donc de montrer le r\'{e}sultat pour $H_{{\frak M}}^{p}(C)$. Pour cela on
distingue les cas suivants :

\begin{enumerate}
\item  $0\leqslant p\leqslant n_1$, o\`{u} $n_1\leqslant n_2\leqslant
...\leqslant n_s$.

Il est clair que $(\sigma _{i}=X_{1,i}...X_{s,i})_{1\leqslant i\leqslant
n_{1}+1}$une $C-$ suite, d'o\`{u} $prof_{{\frak M}}\left( C\right) \geqslant
n_{1}+1$, et d'apr\`{e}s \cite{MTMURA}, on a $H_{{\frak M}}^{p}\left(
C\right) =0,$ pour $0\leqslant p\leqslant n_{1}.$

\item  Lorsque $p>n_{1},$on sait d'apr\`{e}s \ref{EXPL}, du fait que $%
n_{j}\geqslant 1$ pour tout $j, 1\leqslant j\leqslant s$ que 
\begin{equation}
H^{p}(U_{j},O_{U_{j}})=\left\{ 
\begin{array}{ll}
C_{j} & \,\,\,\,\text{si }p=0 \\ 
H_{{\frak M}_{j}}^{n_{j}+1}(C_{j}) &  \,\,\,\,\text{si }p=n_{j}
\\ 
0 & \,\,\,\,\text{si }p\notin \left\{
0,n_{j}\right\} \text{ }
\end{array}
\right.  \label{F*}
\end{equation}
o\`u 
$$H_{{\frak M}_{j}}^{n_{j}+1}(C_{j})  =\frac{1}{X_{j,1}...X_{j,n_{j}+1}}A%
\left[ X_{j,1}^{-1},...,X_{j,n_{j}+1}^{-1}\right]$$
est un $A-$ module libre. Et comme $U=U_{1}\times _{S}...\times _{S}U_{s},$
on a, d'apr\`{e}s la formule de K\"{u}nneth, 
\[
H^{p}(U,O_{U})=\bigoplus\limits_{p_{1}+\ldots
+p_{s}=p}\bigotimes\limits_{j=1}^{s}H^{p_{j}}(U_{j},O_{U_{j}}), 
\]
et $H_{{\frak M}}^{p}\left( C\right) =H^{p-1}(U,O_{U}),$ pour $p\geqslant 2.$
On en d\'{e}duit
\end{enumerate}

\begin{center}
\begin{equation}
H_{{\frak M}}^p\left( C\right) =\bigoplus_{\begin{array}{l} {\small { p_1+...+p_s=p-1}}  \\{\small p_j=0 \text{ ou }n_j } \end{array}} \bigotimes_{j=1}^sH^{p_j}(U_j,O_{U_j}).
\label{For-Kuneth}
\end{equation}
\end{center}

Comme $p>n_{1}$, il r\'{e}sulte de la formule\ref{F*} que $^{^{\prime
}}E^{l,p}=0$, pour tout $p\neq 1+\sum\limits_{j\in J}n_{j}$ o\`{u} $J$ est
une partie non vide de $\{1,...,s\}$ 
\endproof%
%

\begin{coro}
Pour $l+p>N=1+\sum\limits_{j=1}^sn_j=\max\limits_{J\subseteq \left\{
1,\ldots ,s\right\} }\left( 1+\sum\limits_{j\in J}n_j\right) ,$ avec $-r\leq
l\leq 0$ on a $^{\prime }E_1^{l,p}=0$.
\end{coro}

Comme $\underline{f}=(f_{1},...,f_{r})$ est une $C_{\sigma _{i}}-$ suite,
elle est \'{e}galement une $\stackrel{\vee }{{\cal C}^{p}}-$ suite pour $%
p\neq 0$. Or $K^{\bullet ,p}$ est le complexe de Koszul associ\'{e} \`{a} la
suite $\underline{f}$ et au module $\stackrel{\vee }{{\cal C}^{p}};$il est
donc acyclique pour $n\geq 1$ d'apr\`{e}s \cite{MTMURA}, et on a :

\begin{propo}
Pour $p\neq 0,$ on a $H^p(K^{\bullet ,i})=0,$ pour tout $i\in \{1,...,q\}.$
\end{propo}

On d\'{e}duit aussi qu'on a : 
\begin{equation}
^{^{\prime \prime }}E_{2}^{l,p}=\left\{ 
\begin{array}{ll}
0 & \text{si }p\neq 0\text{ et }l\neq 0 \\ 
H_{{\frak M}}^{l}\left( B\right) & \text{si }p=0 \\ 
H^{p}\left( K^{\bullet }\right) & \text{si }p\neq 0\text{ et }l=0.
\end{array}
\right.  \label{Supp-SS2}
\end{equation}

En vertu du r\'{e}sultat ci-dessus et des propri\'{e}t\'{e}s de cette suite
spectrale $\left( ^{^{\prime \prime }}E\right) $ associ\'{e}e \`{a} un
bicomplexe limit\'{e} sup\'{e}rieurement et inf\'{e}rieurement on a $%
^{^{\prime \prime }}E_{2}^{l,p}=$ $^{^{\prime \prime }}E_{3}^{l,p}=\cdots =$ 
$^{^{\prime \prime }}E_{\infty }^{l,p}=E^{l+p}$ o\`{u} $E^{l+p}$ est
l'aboutissement$.$ Par cons\'{e}quent, on a les r\'{e}sultats suivants :

\begin{propo}
\label{About SS}

\begin{enumerate}
\item  $E^m=H^m\left( K^{\bullet }\right) $ si $m<0$ et $E^m=H_{{\frak M}%
}^m\left( B\right) $ si $m\geqslant 0,$

\item  $E^m=0$ si $m>N$.

\item  $H_{{\frak M}}^m\left( B\right) =0$ si $m>N$.
\end{enumerate}
\end{propo}

\subsection{Etudes des groupes de cohomologie $H_{{\frak M}}^{\bullet
}\left( B\right) $ et $H^{\bullet }\left( K^{\bullet }\right) $}

On se propose dans cette partie de g\'{e}n\'{e}raliser le th\'{e}or\`{e}me
de Hurwitz am\'{e}lior\'{e} par J.P.Jouanolou \cite{JONLOU2} dans le cas
d'un seul paquet au cas de $s$ paquets de variables.

\begin{propo}
Supposons que le nombre $r$ de polyn\^{o}mes plurihomog\`{e}nes $f_1,\ldots
,f_r$ est inf\'{e}rieur ou \'{e}gal au nombre de variables du premier paquet
( $r\leqslant n_1+1).$ Alors :

\begin{enumerate}
\item  le complexe de Koszul est acyclique sauf en degr\'{e} $0$ : $%
H^i\left( K^{\bullet }\right) =0$ si $i\neq 0$

\item  si $r<n_1+1$ alors $H_{{\frak M}}^i\left( B\right) =0$ pour $%
i=0,\cdots ,n_1+1-r$.
\end{enumerate}
\end{propo}

{\it \noindent D\'{e}monstration}{\bf . }Des caract\'{e}risations des
supports des deux suites spectrales $^{^{\prime }}E$ et $^{^{\prime \prime
}}E$ (\ref{première pro suite spec} et \ref{Supp-SS2}) on d\'{e}duit que :

\begin{enumerate}
\item  Si $r<n_{1}+1,$ alors, pour tout couple $(l,p)$ d'entiers relatifs
tels que $l+p<n_{1}+1-r,$ on a $^{^{\prime }}E_{1}^{l,p}=0$, par suite $%
E^{i}=0,\,$pour $i<n_{1}+1-r$ , et, donc, $H_{{\frak M}}^{i}\left( B\right)
=0$ pour $0\leqslant i<n_{1}+1-r$ et $H^{i}\left( K^{\bullet }\right) =0$,
pour $0<i<n_{1}+1-r$.

\item  Si $r=n_{1}+1$ alors on a $^{^{\prime }}E_{1}^{l,p}=0$ pour $l+p<0$
et $E^{i}=0$ pour $i<0$.
\end{enumerate}

Ceci montre que $H^{i}\left( K^{\bullet }\right) =0$ si $i\neq 0.$ 
\endproof%
%

\noindent {\bf Remarque} D'apr\`es la d\'efinition 3.1 on a 
$${\cal{T}}=\pi^{-1}(H_{{\frak M}}^0\left( B\right)).$$
\begin{coro}
Si $r<n_1+1$ alors l'id\'{e}al des formes d'inertie est donn\'{e} par ${\cal %
T=}\left( f_1,...,f_r\right) ,$ en particulier, l'id\'{e}al r\'{e}sultant $%
{\frak A}$ est nul. 
\end{coro}

Rappelons (\ref{DN3}) que les mon\^{o}mes $(\sigma _{i_1,\cdots ,i_s})$ sont
plurihomog\`{e}nes de degr\'{e} $\left( 1,\cdots ,1\right) $, ce qui induit
une ${\Bbb Z}^s-$graduation sur les complexes colonnes $K^l\bigotimes%
\nolimits_C\stackrel{\vee }{{\cal C}}^{\bullet }$ du diagramme (\ref{Diag de
comp}) , et, puisque, les diff\'{e}rentielles $d^{\prime \prime }$ sont
plurihomog\`{e}nes de degr\'{e}s $\left( 0,\cdots ,0\right) ,$ alors les
groupes de cohomologie sont \'{e}galement ${\Bbb Z}^s-$gradu\'{e}s.

Pour tous $l=0,\ldots ,r$ et $p=1,\ldots ,s,$ on pose, 
\[
\left\{ 
\begin{array}{lll}
\delta _{j}\left( l\right) & = & {\max }_{1\leq i_{1}<\ldots <i_{l}\leq r} \left( d_{i_{1},j}+\ldots +d_{i_{l},j}\right) -n_{j}-1 \\ 
\delta _{j}\left( 0\right) & = & -n_{j}-1,
\end{array}
\right. 
\]
on a donc, $\delta _{j}\left( 0\right) <\delta _{j}\left( 1\right) <\ldots
<\delta _{j}\left( r\right) =d_{1,j}+\ldots +d_{r,j}-n_{j}-1.$ On note $%
\delta =\left( \delta _{1},\ldots ,\delta _{s}\right) ,$ o\`{u} $\delta
_{j}= $ $\delta _{j}\left( r\right) $ pour tout $1\leq j\leq s$.

\begin{propo}
\label{Part-homo-SS}Pour tous $l=0,\ldots ,r$ et $p=1+n_1,\ldots ,N,$ on a :

\begin{enumerate}
\item  $\left( ^{^{\prime }}E_{1}^{-l,p}\right) _{\left( 0,\cdots ,0\right)
}=0$ si $p<N=1+\sum\limits_{j=1}^{s}n_{j}$,

\item  $\left( ^{^{\prime }}E_{1}^{-l,p}\right) _{\left( \nu _{1},\cdots
,\nu _{s}\right) }=0$ pour $\nu _{j}>\delta _{j}\left( l\right) .$ $\ $
\end{enumerate}
\end{propo}

{\it \noindent D\'{e}monstration.} Il r\'{e}sulte de la proposition \ref
{première pro suite spec} que, pour tout $p\neq 1+\sum\limits_{j\in J}n_{j},$
o\`{u} $J$ est une partie non vide de $\{1,\cdots ,s\}$, la premi\`{e}re
suite spectrale v\'{e}rifie $^{^{\prime }}E_{1}^{-l,p}=0$; on se ram\`{e}ne
donc au cas de $p=1+\sum\limits_{j\in J}n_{j}.$ Les formules $^{^{\prime
}}E_{1}^{-l,p}=H_{{\frak M}}^{p}\left( K^{-l}\right)
=K^{-l}\bigotimes\nolimits_{C}H_{{\frak M}}^{p}(C)$ (\ref{Suite spec2})
\thinspace \ref{For-Kuneth} donnent 
\[
\left( ^{^{\prime }}E_{1}^{0,p}\right) _{\left( 0,\ldots ,0\right) }=\left(
H_{{\frak M}}^{p}(C)\right) _{\left( 0,\ldots ,0\right) }=\bigoplus \limits_{\begin{array}{l}{ p_{1}+...+p_{s}=p-1}\\ {p_{j}=0\text{ ou }n_{j}}
\end{array}} \bigotimes%
\limits_{j=1}^{s}\left( H^{p_{j}}(U_{j},O_{U_{j}})\right) _{0}\text{,} 
\]
dans ces conditions, il existe $j=1,\cdots ,s$ tel que $p_{j}=n_{j},$
or on a $H^{n_{j}}(U_{j},O_{U_{j}})=0\,$ (d'apr\`{e}s \ref{F*}), par
cons\'{e}quent on a $\left( ^{\prime }E_{1}^{0,p}\right) _{\left( 0,\cdots
,0\right) }=0,$ ce qui montre la premi\`{e}re assertion pour $l=0$.

Supposons maintenant que $l\neq 0.$ Pour que la partie plurihomog\`{e}ne de
degr\'{e} $\left( 0,\ldots ,0\right) $ de la premi\`{e}re suite spectrale
soit nulle $(\left( ^{^{\prime }}E_{1}^{-l,p}\right) _{\left( 0,\cdots
,0\right) }=0)$ il faut et il suffit que l'on ait $\left( H_{{\frak M}%
}^{p}\left( C\right) \right) _{-\underline{d}_{i_{1}}-\cdots -\underline{d}%
_{i_{l}}}=0$ pour tous $1\leq i_{1}\leq \cdots \leq i_{l}\leq r$ (cf \ref
{Comp Kosz2}). La relation 
\[
\left( H_{{\frak M}}^{p}\left( C\right) \right) _{-\underline{d}%
_{i_{1}}-\cdots -\underline{d}_{i_{l}}}=\bigoplus\limits_{\begin{array}{l}{ p_{1}+.+p_{s}=p-1}\\ { p_{j}=0\text{ ou }n_{j}}
\end{array}} 
 \left(
\bigotimes\limits_{j=1}^{s}\left( H^{p_{j}}(U_{j},O_{U_{j}})\right)
_{-d_{i_{1},j}-\cdots -d_{i_{l},j}}\right) 
\]
et le fait que $p<N$ montrent qu'il existe $j\in \{1,\cdots ,s\}$ tel que $%
p_{j}=0,$ donc d'apr\`{e}s \ref{F*} $\left( H^{0}(U_{j},O_{U_{j}})\right)
_{-d_{i_{1},j}-\cdots -d_{i_{l},j}}=0\,$ ce qui ach\`{e}ve la
d\'{e}monstration de la premi\`{e}re assertion.

Montrons maintenant la seconde assertion${\bf .}$ Soient $%
p=1+\sum\limits_{j\in J}n_j$ o\`{u} $J$ est une partie non vide de $%
\{1,\cdots ,s\}$ et $\nu =\left( \nu _1,\cdots ,\nu _s\right) \in {\Bbb Z}%
^s. $ Nous avons la formule qui met en relief les cas $l=0$ et $l\neq 0$ :

\[
\left( ^{^{\prime }}E_{1}^{-l,p}\right) _{\nu }=\left\{ 
\begin{array}{ll}
\left( H_{{\frak M}}^{p}\left( C\right) \right) _{\nu } & \text{si }l=0 \\ 
\bigoplus\limits_{1\leqslant i_{1}<...<i_{l}\leqslant r}\left( H_{{\frak M}%
}^{p}\left( C\right) \right) _{\nu -\underline{d}_{i_{1}-\cdots -\underline{d%
}_{i_{l}}}} & \text{si }l\neq 0
\end{array}
\right. \text{.} 
\]

Lorsque $l=0$, la formule \ref{For-Kuneth} entra\^{i}ne

\[
\left( ^{^{\prime }}E_{1}^{0,p}\right) _{\nu }=\bigoplus\limits_{ \begin{array}{l} p_{1}+...+p_{s}=p-1\\  p_{j}=0\text{ ou }n_{j}
\end{array}} \bigotimes\limits_{j=1}^{s}%
\left( H^{p_{j}}(U_{j},O_{U_{j}})\right) _{\nu _{j}}\text{,} 
\]
et puisque $p\geq n_{1}+1,$ il existe $j=1,\cdots ,s$ tel que $p_{j}=n_{j}$,
donc d'apr\`{e}s \ref{EXPL} on a $\left( H^{n_{j}}(U_{j},O_{U_{j}})\right)
_{\nu _{j}}=0$ si $\nu _{j}>-n_{j}-1$, d'o\`{u} $\left( ^{^{\prime
}}E_{1}^{0,p}\right) _{\nu }=0$ si $\nu _{j}>-n_{j}-1$ pour tout $1\leq
j\leq s$.

Lorsque $l\neq 0$, on observe que

\[
\left( ^{^{\prime }}E_{1}^{-l,p}\right) _{\nu }=\bigoplus\limits_{1\leqslant
i_{1}<...<i_{l}\leqslant r}\left( ^{^{\prime }}E_{1}^{0,p}\right) _{\nu -%
\underline{d}_{i_{1}-\cdots -\underline{d}_{i_{l}}}}\text{.} 
\]
Or $\left( ^{^{\prime }}E_{1}^{0,p}\right) _{\nu -\underline{d}%
_{i_{1}}-\cdots -\underline{d}_{i_{l}}}=0$ si $\nu _{j}>d_{i_{1},j}+\cdots
+d_{i_{l},j}-n_{j}-1$ et $j=1,\cdots ,s.$ Par cons\'{e}quent

$\left( ^{\prime }E_1^{-l,p}\right) _\nu =0$, si $\nu _j>{\max }_{1\leq
i_1<\cdots <i_l\leq r} \left( d_{i_1,j}+\cdots +d_{i_l,j}\right)
-n_j-1=\delta _j\left( l\right) $ pour tout $1\leq j\leq s,$ ce qui montre
la proposition. 
\endproof%
%

\begin{propo}
$r$ \'{e}tant le nombre de polyn\^{o}mes plurihomog\`{e}nes $f_1,\ldots ,f_r$%
, on a :

\begin{enumerate}
\item  $\left( H^i\left( K^{\bullet }\right) \right) _{\left( 0,\cdots
,0\right) }=0$ si $r\leq N$ et $i\neq 0$,

\item  $\left( H_{{\frak M}}^i\left( B\right) \right) _{\left( 0,\cdots
,0\right) }=0$ si $r<N$ et $0\leq i<N-r$.
\end{enumerate}
\end{propo}

\noindent {\it D\'{e}monstration. }La proposition pr\'{e}c\'{e}dente donne $%
\left( ^{\prime }E_{1}^{l,p}\right) _{\left( 0,\cdots ,0\right) }=0$ si $%
p<N. $ Ainsi pour tous $l=-r,\cdots ,0$ et $p=1+n_{1},\cdots ,N$ tels que $%
l+p<N-r,$ on a $\left( ^{^{\prime }}E_{1}^{l,p}\right) _{\left( 0,\cdots
,0\right) }=0,$ ce qui montre que l'aboutissement en degr\'{e} $\left(
0,\cdots ,0\right) $ est nul : $E_{\left( 0,\cdots ,0\right) }^{i}=0$ si $%
i<N-r$. La proposition r\'{e}sulte des relations $E^{m}=H^{m}\left(
K^{\bullet }\right) $ si $m<0$ et $E^{m}=H_{{\frak M}}^{m}\left( B\right) $
si $m\geq 0;$ ceci ach\`{e}ve la d\'{e}monstration 
\endproof%
%

\begin{coro}
Sous les hypoth\`{e}ses et notations ci dessus, on a :

\begin{enumerate}
\item  $\left( H_{{\frak M}}^{i}\left( B\right) \right) _{\nu }=0$ pour tout 
$i,$

\item  $\left( H^{i}\left( K^{\bullet }\right) \right) _{\nu }=0$ pour tout $%
i\neq 0.$
\end{enumerate}

pour tout $\nu =\left( \nu _1,\ldots ,\nu _s\right) \in {\Bbb N}^s$ tel que $%
\nu _j>\delta _j$.
\end{coro}

\begin{coro}
Pour $\nu =\left( \nu _1,\ldots ,\nu _s\right) \in {\Bbb N}^s$ avec $\nu
_j>\delta _j,$ on a ${\cal T}_\nu {\cal =}\left( f_1,\ldots ,f_r\right) _\nu
.$
\end{coro}

\subsubsection{Etudes de $\left( H_{{\frak M}}^{\bullet }\left( B\right)
\right) _\protect\protect\delta $ et $\left( H^{\bullet }\left( K^{\bullet
}\right) \right) _\protect\protect\delta $}

L'\'{e}tude des parties plurihomog\`{e}nes $\left( H_{{\frak M}}^{\bullet
}\left( B\right) \right) _\delta $ et $\left( H^{\bullet }\left( K^{\bullet
}\right) \right) _\delta $ des groupes de cohomologies n\'{e}cessite la :

\begin{propo}
On a 
\[
\left( ^{^{\prime }}E_{1}^{-l,p}\right) _{\delta }=\left\{ 
\begin{array}{ll}
A & \text{si }(l,p)=(r,N) \\ 
0 & \text{si }(l,p)\neq (r,N).
\end{array}
\right. 
\]
\end{propo}

\noindent {\it D\'{e}monstration. }On sait que, pour $(l,p)\notin \{0,\ldots
,r\}\times \{1+n_{1},\ldots ,N\},$ $^{^{\prime }}E_{1}^{-l,p}=0.$ Il suffit
de montrer la proposition pour $(l,p)\in \{0,\ldots ,r\}\times
\{1+n_{1},\ldots ,N\}._{{}}$Rappelons que $\delta =\left( \delta _{1},\ldots
,\delta _{s}\right) \in {\Bbb Z}^{s}$ o\`{u} $\delta
_{j}=\sum\limits_{i=1}^{r}d_{i,j}-n_{j}-1.$

Pour $0\leq l<r,$ on a $\delta _{j}\left( l\right) <\delta _{j},$ donc $%
(^{^{\prime }}E_{1}^{-l,p})_{\delta _{j}}=0$ en vertu de la proposition \ref
{Part-homo-SS}, $r$ \'{e}tant le nombre de polyn\^{o}mes plurihomog\`{e}nes $%
f_{1},\ldots ,f_{r}$.

Pour $l=r,$ on sait que $^{^{\prime }}E_{1}^{-r,p}=H_{{\frak M}}^{p}\left(
K^{-r}\right) =H_{{\frak M}}^{p}\left( C\right) [-$\underline{$d$}$%
_{1}-\ldots -$\underline{$d$}$_{r}],$ donc $\left( ^{^{\prime
}}E_{1}^{-r,p}\right) _{\delta }=\left( H_{{\frak M}}^{p}\left( C\right)
\right) _{\left( -n_{1}-1,\ldots ,-n_{s}-1\right) }$ et la formule de
K\"{u}nneth donne 
\[
\left( ^{^{\prime }}E_{1}^{-r,p}\right) _{\delta }=\bigoplus\limits_{\begin{array}{l}
p_{1}+...+p_{s}=p-1  \\ p_{j}=0\text{ ou }n_{j}  \end{array}} \bigotimes%
\limits_{j=1}^{s}H^{p_{j}}(U_{j},O_{U_{j}})_{-n_{j}-1}. 
\]
Si $p<N,$ alors on peut supposer, d'apr\`{e}s la proposition \ref{première
pro suite spec}, que $p=1+\sum\limits_{j\in J}n_{j}$ o\`{u} $J$ est une
partie non vide de $\{1,\cdots ,s\},$ et puisque dans l'expression de $%
\left( ^{^{\prime }}E_{1}^{-r,p}\right) _{\delta }$ ci-dessus $p_{j}$ ne
prend que deux valeurs $0$ ou $n_{j}$ alors il existe $j\in \{1,\ldots ,s\}$
tel que $p_{j}=0.$ Comme $H^{0}\left( U_{j},O_{U_{j}}\right) _{-n_{j}-1}=0$ (%
\ref{F*}), on d\'{e}duit que $\left( ^{^{\prime }}E_{1}^{-r,p}\right)
_{\delta }=0.$

Lorsque $p=N,$ on a $\left( ^{^{\prime }}E_{1}^{-r,N}\right) _{\delta
}=\bigotimes\limits_{j=1}^{s}H^{n_{j}}(U_{j},O_{U_{j}})_{-n_{j}-1}.$ De la
formule \ref{F*}, on d\'{e}duit que $%
H^{n_{j}}(U_{j},O_{U_{j}})_{-n_{j}-1}=A, $ par cons\'{e}quent $\left(
^{^{\prime }}E_{1}^{-r,N}\right) _{\delta }=A.$ 
\endproof%
%

De la proposition pr\'{e}c\'{e}dente on d\'{e}duit l'aboutissement de la
premi\`{e}re suite spectrale en degr\'{e} $\delta $ :

\[
E_\delta ^i=\left\{ 
\begin{array}{ll}
A & \text{si }i=N-r \\ 
0 & \text{si }i\neq N-r.
\end{array}
\right. 
\]
En comparant la formule ci dessus avec 
\[
E_\delta ^i=\left\{ 
\begin{array}{ll}
\left( H_{{\frak M}}^i\left( B\right) \right) _\delta & \text{si }i\geq 0 \\ 
\left( H^i\left( K^{\bullet }\right) \right) _\delta & \text{si }i\leq 0,
\end{array}
\right. 
\]
on d\'{e}duit le th\'{e}or\`{e}me suivant, \thinspace qui est un bilan des
r\'{e}sultats pr\'{e}c\'{e}dents. Rappelons que $r$ est le nombre de
polyn\^{o}mes g\'{e}n\'{e}riques plurihomog\`{e}nes $f_i$ et $N=1+n_1+\ldots
+n_s$ o\`{u} $1+n_j$ est le nombre de variables dans le paquet $j.$

\begin{theo}
Dans les hypoth\`{e}ses et notations ci dessus on a :

\begin{enumerate}
\item  Si $r<N,$ alors

\begin{enumerate}
\item  le complexe $K_\delta ^{\bullet }$ (\ref{Comp Kosz2}) est acyclique
sauf en degr\'{e} $0$ : $\left( H^i(K^{\bullet })\right) _\delta =0$ pour $%
i\neq 0,$

\item  Lorsque $i\geq 0$ et $i\neq N-r,$ on a $\left( H_{{\frak M}}^i\left(
B\right) \right) _\delta =0.$
\end{enumerate}

\item  Pour $r=N$ on a :

\begin{enumerate}
\item  le complexe $K_\delta ^{\bullet }$ est acyclique sauf en degr\'{e} $0$
: $\left( H^i(K^{\bullet })\right) _\delta =0$ pour $i\neq 0,$

\item  le $A-$module $\left( H_{{\frak M}}^0\left( B\right) \right) _\delta $
est libre de rang $1$ et $\left( H_{{\frak M}}^i\left( B\right) \right)
_\delta =0,$ pour $i>0.$
\end{enumerate}

\item  Supposons que $r>N$

\begin{enumerate}
\item  Lorsque $i\neq 0$, on a $\left( H^{i}(K^{\bullet })\right) _{\delta
}=A$ pour $i=N-r,$ et $\left( H^{i}(K^{\bullet })\right) _{\delta }=0$ pour $%
i\neq N-r$

\item  Pour $i\geq 0$, on a $\left( H_{{\frak M}}^i\left( B\right) \right)
_\delta =0.$
\end{enumerate}
\end{enumerate}
\end{theo}

\begin{coro}
On a :

\begin{enumerate}
\item  Si $r\neq N,$ la partie plurihomog\`{e}ne de degr\'{e} $\delta $ de
l'id\'{e}al des formes d'inertie est ${\cal T}_\delta =\left( f_1,\ldots
,f_r\right) _\delta $

\item  Si $r\leq N$ alors le complexe $K_\delta ^{\bullet }$ d\'{e}finit une
r\'{e}solution libre du $A-$module $B_\delta .$
\end{enumerate}
\end{coro}

Dans le cas o\`{u} le nombre de polyn\^{o}mes $f_1,\ldots ,f_r$ est \'{e}gal
\`{a} $N=1+n_1+\ldots ,n_r$ , le $A-$module $\left( H_{{\frak M}}^0\left(
B\right) \right) _\delta $ est libre de rang $1,$ dont on donne ici un
g\'{e}n\'{e}rateur sans donner le lien avec le d\'{e}terminant $"$Jacobien$"$
des polyn\^{o}mes $f_1,\ldots ,f_r,$ \'{e}tabli dans \cite{CHKR}

On d\'{e}finit par r\'{e}currence sur $j=1,\ldots ,s$, des polyn\^{o}mes
uniques $f_{i,l}^{\left( 1\right) }$ (o\`{u} $1\leq l\leq n_{j}+1$) par : 
\begin{equation}
\begin{array}{ll}
f_{i} & =X_{1,1}f_{i,1}^{\left( 1\right) }+\ldots
+X_{1,n_{1}+1}f_{i,n_{1}+1}^{\left( 1\right) } \\ 
f_{i,l}^{\left( 1\right) } & \in A\left[ X_{1,l},\ldots ,X_{1,n_{1}+1}\right]
\left[ \underline{X}_{2},\ldots ,\underline{X}_{s}\right] ,
\end{array}
\label{Décompo-f}
\end{equation}
pour tout $i=1,\ldots ,r=N.$ Pour d\'{e}finir $f_{i,l}^{(j)},$ pour $2\leq
j\leq s$ et $1\leq l\leq n_{j}+1,$ on d\'{e}compose $f_{i,n_{j-1}+1}^{(j-1)}$
sous la forme \ref{Décompo-f} par rapport au paquet \underline{$X$}$_{j}$ : 
\[
\begin{array}{ll}
f_{i,n_{j-1}+1}^{(j-1)} & =X_{j,1}f_{i,1}^{(j)}+\ldots
+X_{j,n_{j}+1}f_{i,n_{j}+1}^{(j)} \\ 
f_{i,l}^{(j)} & \in A\left[ X_{1,n_{1}+1},\ldots ,X_{j-1,n_{j-1}+1}\right] 
\left[ X_{j,l},\ldots ,X_{j,n_{1}+1}\right] \left[ \underline{X}%
_{j+1},\ldots ,\underline{X}_{s}\right] ,
\end{array}
\]
on obtient alors une d\'{e}composition

$f_{i}=\sum\limits_{l=1}^{n_{1}}X_{1,l}f_{i,l}^{\left( 1\right)
}+X_{1,n_{1}+1}\sum\limits_{l=1}^{n_{2}}f_{i,l}^{\left( 2\right) }+\ldots
+\prod\limits_{k=1}^{j-1}X_{k,n_{k}+1}\sum%
\limits_{l=1}^{n_{j}}x_{j,l}f_{i,l}^{(j)}+\ldots
+\prod\limits_{k=1}^{s-1}X_{k,n_{k}+1}\sum%
\limits_{l=1}^{n_{s}}x_{s,l}f_{i,l}^{(s)}+\prod%
\limits_{k=1}^{s}X_{k,n_{k}+1}f_{i,n_{s}+1}^{(s)}$ du polyn\^{o}me $f_{i}$
pour $i=1,\ldots ,r.$ Consid\'{e}rons le d\'{e}terminant d'ordre $N$ 
\[
{\cal D=}\left| 
\begin{array}{lllllllllll}
f_{1,1}^{\left( 1\right) } & \ldots & f_{1,n_{1+}1}^{\left( 1\right) } & 
\ldots & f_{1,1}^{(j)} & \ldots & f_{1,n_{j}}^{(j)} & \ldots & 
f_{1,1}^{\left( s\right) } & \ldots & f_{1,n_{s}+1}^{\left( s\right) } \\ 
\vdots &  & \vdots &  & \vdots &  & \vdots &  & \vdots &  & \vdots \\ 
f_{i,1}^{\left( 1\right) } & \ldots & f_{i,n_{1}}^{\left( 1\right) } & \ldots
& f_{i,1}^{(j)} & \ldots & f_{i,n_{j}}^{(j)} & \ldots & f_{i,1}^{\left(
s\right) } & \ldots & f_{i,n_{s}+1}^{\left( s\right) } \\ 
\vdots &  & \vdots &  & \vdots &  & \vdots &  & \vdots &  & \vdots \\ 
f_{N,1}^{\left( 1\right) } & \ldots & f_{N,n_{1}}^{\left( 1\right) } & \ldots
& f_{N,1}^{(j)} & \ldots & f_{N,n_{j}}^{(j)} & \ldots & f_{N,1}^{\left(
s\right) } & \ldots & f_{N,n_{s}+1}^{\left( s\right) }
\end{array}
\right| . 
\]
Si on pose $N_{j}=\sum\limits_{l=j+1}^{s}n_{l}$ pour $1\leq j\leq s-1$ et $%
N_{s}=0$ alors on a :

\begin{theo}
\begin{enumerate}
\item  Le d\'{e}terminant ${\cal D}$ est une forme d'inertie
plurihomog\`{e}ne de degr\'{e} $\delta _j-N_j$ par rapport au paquet 
\underline{$X$}$_j,$ c'est \`{a} dire ${\cal D\in T}_{\left( \delta
_1-N_1,\ldots ,\delta _s-N_s\right) }.$

\item  La classe $\Delta $ de $X_{1,n_1+1}^{N_1}\ldots
X_{1,n_{s-1}+1}^{N_{s-1}}{\cal D}$ modulo $\left( f_1,\ldots ,f_N\right) $
est un g\'{e}n\'{e}rateur de $\left( H_{{\frak M}}^0\left( B\right) \right)
_\delta .$
\end{enumerate}
\end{theo}

{\em \ }

\end{document}